\documentclass[final,3p,times]{elsarticle}
\usepackage{amssymb}
\usepackage{amsthm}
\usepackage{mathrsfs}
\usepackage{graphicx}
\usepackage[colorinlistoftodos]{todonotes}
\usepackage[colorlinks=true, allcolors=blue]{hyperref}
\usepackage{amsmath}
\usepackage{multirow}
\usepackage{amsfonts, geometry}
\usepackage{bbold}
\usepackage{subcaption}
\usepackage{booktabs}
\usepackage{url}
\usepackage{tikz}
\usetikzlibrary{positioning}
\usepackage{makecell}
\usepackage{tabularx}
\usepackage{booktabs}

\theoremstyle{dgdef}

\let\today\relax
\makeatletter
\def\ps@pprintTitle{%
\let\@oddhead\@empty
\let\@evenhead\@empty
\def\@oddfoot{}%
\let\@evenfoot\@oddfoot}
\makeatother
\begin{document}
\begin{frontmatter}
\title{Solving the Fisher nonlinear differential equations via Physics-Informed Neural Networks: A Comprehensive Retraining Study and Comparative Analysis with the Finite Difference Method}
\author[1]{Ahmed Aberqi}
\ead{ahmed.aberqi@usmba.ac.ma}
\author[2]{Ahmed Miloudi}
\ead{ahmed.miloudi@usmba.ac.ma}
\address[1]{Laboratory LSATE, National School of Applied Sciences, Sidi Mohamed Ben Abdellah University, Fez, Morocco.}

\address[2]{ Faculté de médecine et de pharmacie, Sidi Mohamed Ben Abdellah University, Fez, Morocco.}

\cortext[cor1]{ Corresponding authors: Ahmed Aberqi (ahmed.aberqi@usmba.ac.ma) and Ahmed Miloudi (ahmed.miloudi@usmba.ac.ma)}

\begin{abstract}
Physics-Informed Neural Networks (PINNs) represent a groundbreaking paradigm in scientific computing, seamlessly integrating the robust framework of deep learning with fundamental physical laws. This paper meticulously applies the standard PINN framework to solve the challenging one-dimensional nonlinear Fisher-KPP equation, a critical model in reaction-diffusion dynamics describing phenomena such as population spread and flame propagation. We detail a comprehensive methodology, encompassing the neural network architecture, the physics-informed loss function, and an in-depth investigation into retraining strategies aimed at optimizing model performance. Our approach is rigorously validated through a direct comparison of the PINN solution against both the known analytical solution and a numerical solution derived from the Finite Difference Method (FDM). Through this work, we elucidate the intricate balance between model complexity, training efficiency, and accuracy. Results highlight the PINN's remarkable capability in accurately approximating the solution to this complex PDE, while also shedding light on the critical aspects and challenges of model retraining, particularly concerning the optimizer's state. This study provides a thorough quantitative error analysis, demonstrating the efficacy of PINNs as a viable and competitive alternative to traditional numerical methods for solving nonlinear differential equations, and discusses their broader applications across various scientific domains.
\end{abstract}


\date{\today}
\end{frontmatter}

\section{Introduction}
The resolution of differential equations, encompassing both ordinary differential equations (ODEs) and partial differential equations (PDEs), forms the bedrock of scientific and engineering endeavors, underpinning simulations and theoretical advancements across disciplines from fluid dynamics and quantum mechanics to finance and epidemiology \cite{A1, A2, A3}. Traditional numerical methodologies, such as the Finite Difference Method (FDM) \cite{Leveque2007}, Finite Element Method (FEM), and Finite Volume Method (FVM), have historically provided robust and reliable frameworks for approximating solutions. These methods rely on discretizing the computational domain into a mesh, which can become computationally prohibitive and geometrically complex in high-dimensional spaces or for problems involving intricate boundaries, a challenge often referred to as the curse of dimensionality \cite{han2018solving, A4}. Furthermore, their explicit reliance on mesh generation can introduce significant overhead and limitations in adaptability to evolving domains or highly irregular geometries.

The recent renaissance in deep learning has ushered in a transformative era for scientific computing \cite{karniadakis2021physics}. Deep Neural Networks (DNNs), with their inherent capacity for universal function approximation \cite{hornik1991universal}, offer a compelling alternative to traditional grid-based solvers. Among the most impactful innovations in this domain are Physics-Informed Neural Networks (PINNs), first formally introduced by Raissi, Perdikaris, and Karniadakis in 2019 \cite{Raissi2019}. PINNs represent a novel paradigm that seamlessly integrates the powerful capabilities of deep learning with fundamental physical principles, typically expressed as PDEs. Unlike purely data-driven machine learning models, PINNs embed the governing equations directly into their loss functions, compelling the neural network to learn solutions that are not only consistent with observed data (if available) but also rigorously adhere to the underlying physical laws. This unique blend renders PINNs highly versatile, enabling them to tackle both forward problems (solving PDEs) and inverse problems (discovering unknown parameters or operators from data) with remarkable efficacy, often in a mesh-free environment.

This paper is dedicated to the comprehensive application of the standard PINN framework to solve the one-dimensional nonlinear Fisher equation, more commonly known as the Fisher-KPP (Kolmogorov–Petrovsky–Piskunov) equation. The Fisher-KPP equation is a canonical reaction-diffusion PDE that models a diverse range of natural phenomena, including population dynamics, the spread of advantageous genes, chemical reaction fronts, and flame propagation \cite{fisher1937wave, kolmogorov1937}. Its intrinsic nonlinearity and the existence of a well-known analytical traveling wave solution make it an invaluable benchmark for evaluating the accuracy, stability, and efficiency of novel numerical and data-driven solution methodologies.

Our study offers several significant contributions. Firstly, we present a detailed implementation and application of a robust PINN architecture tailored for the Fisher-KPP equation under specific initial and boundary conditions. A particular emphasis is placed on elucidating how the governing PDE and its associated constraints are systematically encoded into the PINN's composite loss function, leveraging the power of automatic differentiation. Secondly, we conduct an in-depth experimental investigation into the efficacy and practical implications of systematic retraining strategies on the PINN's performance. This analysis critically examines the nuances and challenges encountered when fine-tuning pre-trained models within a scientific computing context, providing insights into convergence behavior and optimization dynamics. Lastly, we undertake a rigorous comparative analysis, juxtaposing the PINN's derived solution against both the exact analytical solution and a high-fidelity numerical solution obtained via the well-established Finite Difference Method (FDM). This multi-faceted comparison quantitatively assesses the PINN's accuracy, stability, and computational characteristics, positioning it against a traditional solver. Through this work, we aim to underscore the immense potential of PINNs as a powerful, mesh-free, and versatile alternative for solving nonlinear PDEs, while also shedding light on the practical considerations and challenges that accompany their deployment.

\section{Preliminary knowledge: Physics-Informed Neural Networks}

\subsection{Deep Neural Networks as Universal Function Approximators}

At the core of Physics-Informed Neural Networks lies the capability of Deep Neural Networks (DNNs) to approximate complex, nonlinear functions. A DNN is a computational model inspired by the structure and function of the human brain. It consists of multiple layers of interconnected nodes, or "neurons," organized into an input layer, several hidden layers, and an output layer. Each connection between neurons has an associated weight, and each neuron has a bias. The output of a neuron is typically computed by taking a weighted sum of its inputs, adding a bias, and then passing this result through a nonlinear activation function (e.g., sigmoid, ReLU, Tanh).

The theoretical underpinning for using neural networks as function approximators is the Universal Approximation Theorem. Informally, this theorem states that a feedforward neural network with a single hidden layer containing a finite number of neurons (and a non-polynomial activation function) can approximate any continuous function on a compact subset of $\mathbb{R}^n$ to any desired degree of accuracy \cite{hornik1991universal}. Deeper networks (with multiple hidden layers), while not strictly necessary for universal approximation, are empirically found to be more efficient at learning complex, hierarchical representations and can achieve better generalization with fewer parameters for certain tasks. In the context of solving PDEs, a DNN effectively learns a continuous mapping from the spatio-temporal coordinates $(x, t)$ to the solution $u(x, t)$.

The trainable parameters of a neural network are its weights and biases, collectively denoted by $\boldsymbol{\theta}$. The process of "training" a neural network involves iteratively adjusting these parameters to minimize a predefined loss function. This minimization is typically performed using gradient-based optimization algorithms, which rely on computing the gradients of the loss function with respect to the network's parameters.

\subsection{Automatic Differentiation: The Backbone of PINNs}

A critical enabler for PINNs is Automatic Differentiation (AD). Unlike symbolic differentiation (which can be computationally expensive and difficult for complex expressions) or numerical differentiation (which suffers from truncation errors and is sensitive to step size), AD provides exact derivatives by systematically applying the chain rule to the elementary operations that constitute a function. In the context of neural networks, AD propagates gradients through the computational graph of the network, efficiently calculating the derivatives of the network's output (and consequently the loss function) with respect to its input variables and, crucially, with respect to its internal parameters.

For PINNs, AD is indispensable for two main reasons:
\begin{enumerate}
 \item \textbf{PDE Residual Computation:} To enforce the PDE, we need to compute partial derivatives of the network's output $u(x, t)$ with respect to its inputs, such as $\frac{\partial u}{\partial t}$ and $\frac{\partial^2 u}{\partial x^2}$. AD allows these derivatives to be computed precisely and efficiently, forming the residual of the PDE.
 \item \textbf{Parameter Optimization:} The overall loss function of a PINN depends on the network's parameters. AD enables the computation of the gradients $\nabla_{\boldsymbol{\theta}} L$, which are then used by optimizers (like Adam or L-BFGS-B) to update $\boldsymbol{\theta}$ in the direction that minimizes the loss.
\end{enumerate}
Modern deep learning frameworks (e.g., TensorFlow, PyTorch) inherently support AD, making the implementation of PINNs significantly more straightforward and computationally efficient.

\subsection{PINN Architecture and Loss Function}

A Physics-Informed Neural Network $\mathcal{N}(\mathbf{x}; \boldsymbol{\theta})$ serves as a continuous, differentiable approximation of the PDE's solution $u(\mathbf{x})$, where $\mathbf{x}$ represents the input variables (e.g., time $t$ and spatial coordinates $x_1, \dots, x_d$), and $\boldsymbol{\theta}$ denotes the trainable parameters (weights and biases) of the network.

The central tenet of PINNs is encapsulated in their composite loss function, which strategically combines multiple terms to enforce both data fidelity and adherence to physical laws. For a generic PDE $\mathcal{L}u = f$ defined on a domain $\Omega$ with specified initial conditions $u(\mathbf{x},0)=u_0(\mathbf{x})$ and boundary conditions $u(\mathbf{x}_b, t)=u_{BC}(\mathbf{x}_b,t)$ on the boundary $\partial\Omega$, the total loss function $L$ is typically formulated as a weighted sum of three primary components:

\begin{equation}
L = w_{IC} L_{IC} + w_{BC} L_{BC} + w_{Res} L_{Res}
\end{equation}
where $w_{IC}, w_{BC}, w_{Res}$ are weighting coefficients (often initialized to 1, or adaptively adjusted during training to balance the contributions of each term, as discussed in advanced PINN variants):
\begin{itemize}
 \item $L_{IC}$ (Initial Condition Loss): This term quantifies the discrepancy between the network's prediction at the initial time $t=0$ and the true initial state of the system. It is computed as the mean squared error (MSE) over a set of initial condition points $\left\{\mathbf{x}_{i,IC}\right\}_{i=1}^{N_{IC}}$ sampled from the domain at $t=0$. This ensures that the learned solution starts from the correct state. For a 1D problem $u(x,0)=u_0(x)$:
 \begin{equation}
     L_{IC} = \frac{1}{N_{IC}} \sum_{i=1}^{N_{IC}} (u_{NN}(x_{i,IC}, 0) - u_0(x_{i,IC}))^2
\end{equation}
\item $L_{BC}$ (Boundary Condition Loss): This term ensures that the network's solution respects the constraints imposed at the boundaries of the spatial domain. It is calculated as the MSE between the network's prediction at a set of boundary points $\left\{\mathbf{x}_{i,BC}\right\}_{i=1}^{N_{BC}}$ and the prescribed boundary values. This is crucial for defining the problem's spatial confinement. For a 1D problem with Dirichlet boundary conditions at $x_{min}$ and $x_{max}$:
    \begin{equation}
    L_{BC} = \frac{1}{N_{BC}} \sum_{i=1}^{N_{BC}} (u_{NN}(x_{i,BC}, t_{i,BC}) - u_{BC}(x_{i,BC}, t_{i,BC}))^2
    \end{equation}
    \item $L_{Res}$ (Residual Loss): This is the core physics-enforcing term. It represents the mean squared error of the PDE residual, obtained by substituting the network's output $u_{NN}(\mathbf{x})$ into the governing PDE. This term is evaluated at a set of collocation points $\left\{\mathbf{x}_{i,Res}\right\}_{i=1}^{N_{Res}}$ sampled randomly or strategically within the spatio-temporal domain $\Omega$. Minimizing $L_{Res}$ forces the network to satisfy the physical laws throughout the domain.
    \begin{equation}
    L_{Res} = \frac{1}{N_{Res}} \sum_{i=1}^{N_{Res}} (\mathcal{L}u_{NN}(\mathbf{x}_{i,Res}) - f(\mathbf{x}_{i,Res}))^2
 \end{equation}
\end{itemize}
The total loss function is then minimized using gradient-based optimization algorithms.

\subsection{Optimization Strategies and Hyperparameter Tuning}

The optimization of PINNs, like other deep learning models, relies heavily on efficient optimization algorithms and careful hyperparameter tuning. The most commonly used optimizers include:
\begin{itemize}
 \item \textbf{Adam (Adaptive Moment Estimation):} A first-order optimization algorithm that uses adaptive learning rates for each parameter. It computes individual adaptive learning rates for different parameters from estimates of first and second moments of the gradients. Adam is widely popular due to its computational efficiency and good performance across a wide range of problems, often requiring less memory than other adaptive methods \cite{kingma2014adam}.
 \item \textbf{L-BFGS-B (Limited-memory Broyden–Fletcher–Goldfarb–Shanno algorithm with box constraints):} A quasi-Newton optimization algorithm that approximates the Hessian matrix to guide the search for the minimum. L-BFGS-B is a second-order method that can converge faster and often achieve higher accuracy than first-order methods like Adam, especially for problems with smooth loss landscapes. However, it typically requires storing more history of gradients and is not as easily parallelizable, making it more suitable for smaller batch sizes or full-batch optimization \cite{byrd1995limited}.
\end{itemize}
Hyperparameters, such as the learning rate, number of training iterations (epochs), batch size, and the architecture of the neural network (number of layers, neurons per layer, activation functions), critically influence the training process and the final accuracy of the PINN. Successfully training PINNs can be challenging, with potential issues such as slow convergence or failure to train, often due to complex loss landscapes or unbalanced gradient propagation \cite{krishnapriyan2021characterizing}.
\begin{itemize}
\item \textbf{Learning Rate:} This parameter determines the step size at each iteration while moving towards a minimum of the loss function. A small learning rate can lead to slow convergence, while a large one can cause oscillations or divergence. Learning rate schedules (e.g., exponential decay, cosine annealing) are often employed to dynamically adjust the learning rate during training, typically decreasing it over time to allow for finer convergence.
 \item \textbf{Number of Iterations:} The total number of optimization steps. This needs to be balanced against computational cost and the risk of overfitting (though less common in PINNs than purely data-driven models).
\item \textbf{Batch Size:} The number of samples (collocation points, initial/boundary points) used in each iteration to compute the loss and its gradients. Larger batch sizes provide more stable gradient estimates but require more memory and computation per step.
\item \textbf{Loss Weights ($w_{IC}, w_{BC}, w_{Res}$):} Balancing the contributions of different loss components is crucial. If one component dominates, the network might satisfy that condition perfectly but neglect others. While often set to 1, adaptive weighting schemes have been proposed to dynamically adjust these weights during training, for instance, based on gradient magnitudes or magnitudes of the loss terms themselves \cite{wang2021understanding}.
\end{itemize}
Effective hyperparameter tuning often involves a combination of grid search, random search, and empirical judgment based on prior experience with similar problems.

\section{Fisher-KPP Equation: Mathematical Background and Solution Approaches}

\subsection{Mathematical Formulation of the Fisher-KPP Equation}

The one-dimensional Fisher-KPP equation is a fundamental model in mathematical biology, ecology, and chemistry, describing the interplay between diffusive transport and nonlinear reaction. It belongs to a class of reaction-diffusion equations and is given by:

\begin{equation} \label{eq:fisher-kpp}
\frac{\partial u}{\partial t} = D \frac{\partial^2 u}{\partial x^2} + R u (1-u)
\end{equation}
where $u(x,t)$ represents a concentration or population density at position $x$ and time $t$. The equation is defined over a spatio-temporal domain $\Omega = [0, L] \times [0, T]$.
\begin{itemize}
 \item $D > 0$: The diffusion coefficient, which quantifies the rate at which $u$ spreads spatially. A larger $D$ indicates faster diffusion.
\item $R > 0$: The reaction rate, representing the intrinsic growth rate of the population or the rate of reaction. The nonlinear term $R u (1-u)$ is a logistic growth term, meaning that the growth rate is maximal at $u=0.5$ and approaches zero as $u$ approaches 0 or 1. This term models limiting factors such as carrying capacity.
\end{itemize}
For the purpose of this study, we adopted standard parameters commonly found in the literature: $D=0.01$ and $R=1.0$. These fixed values for $D$ and $R$ were directly embedded into the residual loss term of the PINN, ensuring the network learned the solution for this specific instance of the Fisher-KPP equation.

The Fisher-KPP equation is renowned for its traveling wave solutions. These solutions represent patterns that propagate through the domain at a constant speed without changing their shape. For specific initial conditions, the exact analytical solution for a traveling wave is known, providing an invaluable ground truth for validating numerical and approximate methods. For the parameters chosen ($D=0.01, R=1.0$), the exact analytical solution is given by:

\begin{equation} \label{eq:exact-solution}
u(x,t) = \frac{1}{1 + e^{\sqrt{\frac{R}{2D}}(x - \sqrt{2DR}t)}}
\end{equation}
This solution describes a wavefront propagating to the right with a constant speed $c = \sqrt{2DR}$, transitioning smoothly from $u=1$ (or a high concentration/density) to $u=0$ (or a low concentration/density). This specific form of solution is critical for understanding invasion processes and pattern formation in various systems.

\subsection{Numerical Approach: Finite Difference Method (FDM)}

To provide a robust and transparent benchmark for comparison with the PINN solution, we implemented a high-fidelity explicit Finite Difference Method (FDM) to solve the Fisher-KPP equation. FDM is a classical numerical technique that approximates derivatives using finite differences, transforming the continuous PDE into a system of algebraic equations solvable on a discrete grid. The FDM solution serves as a well-understood and reliable baseline against which the performance of the data-driven PINN approach can be critically assessed.

\subsubsection{Numerical Scheme and Stability Analysis}

The FDM scheme discretizes the spatio-temporal domain $x \in [0, 1]$ and $t \in [0, 1]$. We used a uniform spatial grid with $N_x$ points, leading to a spatial step size $\Delta x = L/(N_x-1)$. For this study, we set $N_x=201$, resulting in $\Delta x = 1/(201-1) = 0.005$. For the temporal domain, a uniform grid with $N_t$ steps defines the time step $\Delta t = T/N_t$.

For the Fisher-KPP equation, we employed a standard explicit forward Euler scheme for time advancement coupled with central differences for the spatial second derivative. The time derivative $\frac{\partial u}{\partial t}$ at point $(x_i, t_n)$ is approximated by a forward difference:
$$\frac{\partial u}{\partial t} \approx \frac{u_i^{n+1} - u_i^n}{\Delta t}$$The second spatial derivative $\frac{\partial^2 u}{\partial x^2}$ at point $(x_i, t_n)$ is approximated by a central difference:$$\frac{\partial^2 u}{\partial x^2} \approx \frac{u_{i+1}^n - 2u_i^n + u_{i-1}^n}{(\Delta x)^2}$$
Substituting these approximations into Equation \ref{eq:fisher-kpp} yields the following explicit discrete approximation:

\begin{equation} \label{eq:fdm-scheme}
u_i^{n+1} = u_i^n + \Delta t \left( D \frac{u_{i+1}^n - 2u_i^n + u_{i-1}^n}{(\Delta x)^2} + R u_i^n (1-u_i^n) \right)
\end{equation}
where $u_i^n$ denotes the approximation of $u(x_i, t_n)$ at spatial point $x_i$ and time $t_n$.

Initial conditions $u(x,0)$ and Dirichlet boundary conditions $u(0,t)$ and $u(1,t)$ for the FDM simulation were precisely derived from the exact analytical solution (Equation \ref{eq:exact-solution}) to ensure consistency across comparisons. These conditions are directly applied at the grid points for $t=0$ and at $x=0, x=1$ for all time steps.

A critical aspect of explicit FDM schemes for parabolic PDEs is numerical stability. The time step $\Delta t$ must be chosen carefully to prevent spurious oscillations and ensure convergence. For the diffusion term, the stability criterion is commonly known as the CFL (Courant-Friedrichs-Lewy) condition for diffusion: $\Delta t \le \frac{(\Delta x)^2}{2D}$. The reaction term $R u (1-u)$ can also introduce instability if not handled carefully, potentially requiring even smaller time steps or implicit methods. To ensure robust stability for our explicit scheme while allowing for a reasonable simulation time, we adopted a conservative time step choice. Given $D=0.01$ and $\Delta x=0.005$, the diffusion stability limit is $\Delta t \le \frac{(0.005)^2}{2 \times 0.01} = \frac{0.000025}{0.02} = 0.00125$. We chose $\Delta t = 0.000625$ (i.e., half of the diffusion stability limit), which translates to $N_t = 1.0 / 0.000625 = 1600$ time steps to compute the solution at the final time $t=1.0$. This ensures the numerical solution remains stable and accurate throughout the simulation.

\subsubsection{Error Quantification for Numerical Methods}

To rigorously quantify the performance of both the FDM and PINN solutions, we utilize the relative $L_2$ error. This metric provides a normalized measure of the difference between an approximate solution and the exact analytical solution, allowing for meaningful comparison across different methods and scales. The relative $L_2$ error is defined as:

\begin{equation} \label{eq:l2-error}
E_{L_2,rel} = \frac{||u_{approx} - u_{exact}||_2}{||u_{exact}||_2} = \frac{\sqrt{\sum_{k=1}^{M} (u_{approx,k} - u_{exact,k})^2}}{\sqrt{\sum_{k=1}^{M} (u_{exact,k})^2}}
\end{equation}
where $M$ is the total number of evaluation points across the domain, $u_{approx}$ represents the approximated solution (either FDM or PINN), and $u_{exact}$ is the corresponding analytical solution evaluated at the same points. This metric allows for a direct comparison of the overall accuracy of the numerical solutions against the true solution.

Using this robust FDM scheme, the relative $L_2$ error between the exact analytical solution and the FDM solution at $t=1.0$ was calculated to be $1.42 \times 10^{-4}$. This error serves as a precise and reliable benchmark for the PINN's performance, indicating the inherent accuracy achievable by a well-tuned classical numerical method for this specific problem setup.

\subsection{PINN Approach and Training Configurations}

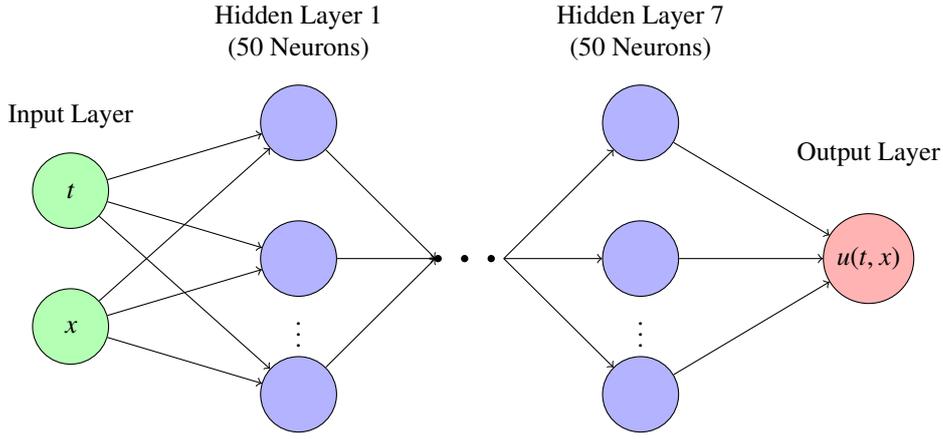
\begin{figure}[h!]
\centering
\begin{tikzpicture}[
 x=1.5cm, y=1.2cm,
neuron/.style={circle, draw, minimum size=1cm, fill=gray!20},
 input/.style={neuron, fill=green!30},
 output/.style={neuron, fill=red!30},
 hidden/.style={neuron, fill=blue!30},
layer_label/.style={align=center}
]
\node[input] (I-1) at (0, 0) {$t$};
 \node[input] (I-2) at (0, -1.5) {$x$};
 \node[layer_label, above=0.2cm of I-1] {Input Layer};

\node[hidden] (H1-1) at (2, 0.75) {};
\node[hidden] (H1-2) at (2, -0.75) {};
\node at (2, -1.5) {\vdots};
\node[hidden] (H1-3) at (2, -2.25) {};
\node[layer_label, above=0.2cm of H1-1] {Hidden Layer 1 \\ (50 Neurons)};

 \node at (3.5, -0.75) {\Huge \dots};
 \node[hidden] (H7-1) at (5, 0.75) {};
\node[hidden] (H7-2) at (5, -0.75) {};
 \node at (5, -1.5) {\vdots};
 \node[hidden] (H7-3) at (5, -2.25) {};
 \node[layer_label, above=0.2cm of H7-1] {Hidden Layer 7 \\ (50 Neurons)};

\node[output] (O-1) at (7, -0.75) {$u(t,x)$};
\node[layer_label, above=0.5cm of O-1] {Output Layer};

\foreach \i in {1,2}
\foreach \j in {1,2,3}
\draw[->] (I-\i) -- (H1-\j);

\foreach \i in {1,2,3}
 \draw[->] (H1-\i) -- (3.2, -0.75);

 \foreach \i in {1,2,3}
 \draw[->] (3.8, -0.75) -- (H7-\i);

 \foreach \i in {1,2,3}
 \draw[->] (H7-\i) -- (O-1);
\end{tikzpicture}
\caption{Schematic of the 7x50 Physics-Informed Neural Network (PINN) architecture. The network takes spatio-temporal coordinates $(t, x)$ as input, processes them through 7 hidden layers of 50 neurons each with Tanh activation functions, and outputs the predicted solution $u(t,x)$.}
\label{fig:pinn_architecture}
\end{figure}

Our PINN model for approximating the solution $u(x,t)$ to the Fisher-KPP equation leverages a deep, fully connected neural network (DNN). The architecture was meticulously designed to balance expressivity with computational efficiency, aiming to capture the complex nonlinear dynamics of the traveling wave solution. Specifically, the network consists of an input layer for the two spatio-temporal coordinates $(t, x)$, followed by \textbf{7 hidden layers, each comprising 50 neurons}, and a single output neuron for $u(x,t)$. This architecture, totaling 7 hidden layers, was chosen after preliminary experiments, providing sufficient depth to learn the complex function mapping without excessive computational burden. The structure of this deep neural network is depicted in Figure \ref{fig:pinn_architecture}.

The hyperbolic tangent (Tanh) activation function, chosen for its differentiability and its ability to introduce non-linearity crucial for approximating complex PDE solutions, was consistently applied between all layers. Tanh is often preferred over ReLU in PINNs because its derivatives are smooth and well-behaved, which is important for the automatic differentiation process used to compute the PDE residual. However, research into adaptive activation functions has shown promise in accelerating convergence and improving accuracy in PINNs \cite{jagtap2020adaptive}. Network weights were initialized using Xavier normal initialization \cite{glorot2010understanding}, a technique known to help maintain signal propagation (preventing vanishing/exploding gradients) through deep networks during training, and biases were initialized to zeros.

The core of our PINN implementation lies in how the governing PDE and its associated conditions are integrated into the loss function. The partial derivatives $\frac{\partial u}{\partial t}$ and $\frac{\partial^2 u}{\partial x^2}$ required by the Fisher-KPP equation (Equation \ref{eq:fisher-kpp}) are automatically computed from the neural network's output $u_{NN}(x,t)$ using TensorFlow's automatic differentiation capabilities. This means that the residual term $L_{Res}$ directly enforces the PDE:
$$\mathcal{R}(u_{NN}) = \frac{\partial u_{NN}}{\partial t} - D \frac{\partial^2 u_{NN}}{\partial x^2} - R u_{NN} (1-u_{NN})$$
The residual loss $L_{Res}$ is then computed as the mean squared error of $\mathcal{R}(u_{NN})$ evaluated at a set of randomly sampled collocation points within the spatio-temporal domain $[0,1] \times [0,1]$. For this study, we used 10,000 randomly sampled collocation points for each training iteration, ensuring a diverse representation of the domain.

Similarly, the initial condition $u(x,0) = \frac{1}{1 + e^{\sqrt{\frac{R}{2D}}(x - 0)}}$ and the Dirichlet boundary conditions $u(0,t) = \frac{1}{1 + e^{\sqrt{\frac{R}{2D}}(- \sqrt{2DR}t)}}$ and $u(1,t) = \frac{1}{1 + e^{\sqrt{\frac{R}{2D}}(1 - \sqrt{2DR}t)}}$ are enforced through their respective MSE loss terms ($L_{IC}$ and $L_{BC}$) by sampling points along the initial time slice ($t=0$) and at the spatial boundaries ($x=0, x=1$). For the initial conditions, 1,000 points were sampled across the spatial domain at $t=0$. For the boundary conditions, 1,000 points were sampled along each spatial boundary ($x=0$ and $x=1$) across the temporal domain.

The training objective for our PINN was to minimize this composite loss function. This multi-objective optimization problem was tackled using the Adam optimizer, renowned for its adaptive learning rate capabilities and computational efficiency in deep learning applications \cite{kingma2014adam}. The initial learning rate for the primary training phase was set to $1 \times 10^{-3}$. An exponential learning rate scheduler, with a decay factor of $0.99$, was also employed to gradually reduce the learning rate over iterations. This scheduling strategy is crucial for stable convergence, allowing the optimizer to make larger steps early in training and then finer adjustments as it approaches the minimum of the loss landscape, preventing oscillations.

\begin{figure}[h!]
    \centering
    \begin{subfigure}[b]{0.49\textwidth}
        \centering
        \includegraphics[width=\textwidth]{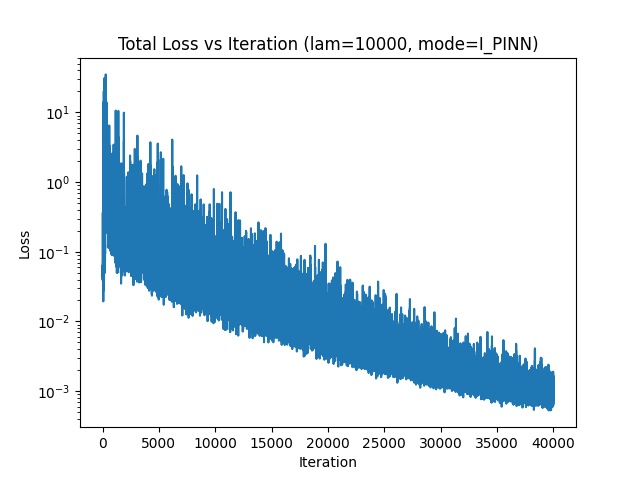}
        \caption{Total Loss vs. Iteration.}
        \label{fig:total_loss}
    \end{subfigure}
    \hfill
    \begin{subfigure}[b]{0.49\textwidth}
        \centering
        \includegraphics[width=\textwidth]{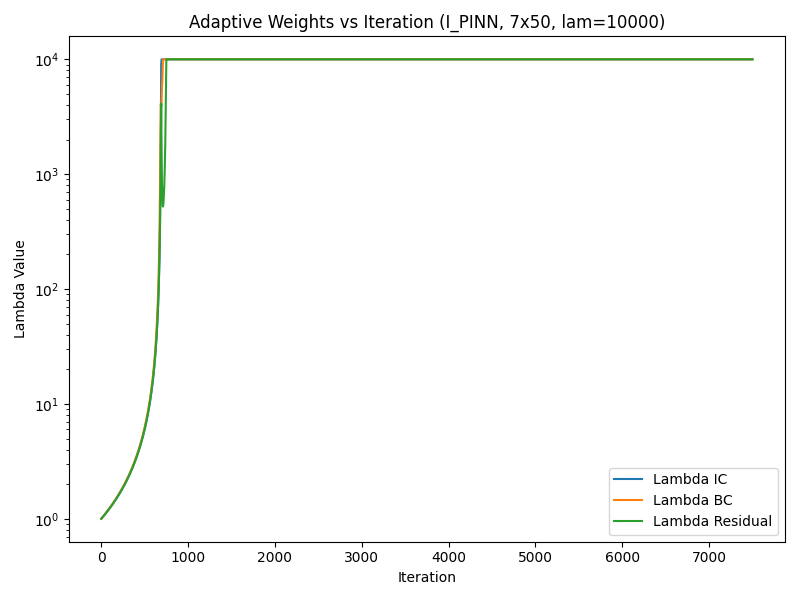}
        \caption{Adaptive Weights (Lambda) vs. Iteration.}
        \label{fig:lambda_plot}
    \end{subfigure}
    \caption{Training dynamics of the I-PINN model with a 7x50 architecture. (a) The total loss function decreases steadily over 40,000 iterations, shown on a logarithmic scale. (b) The adaptive weights for the initial condition (IC) and boundary conditions (BC) quickly saturate at the predefined maximum value of $10^4$, ensuring these conditions are strictly enforced.}
 \label{fig:loss_and_weights}
\end{figure}

Our study involved different training configurations, detailed in Table \ref{tab:training_configs}, designed to explore the performance characteristics of the PINN under various training conditions. For this paper, we also utilized an adaptive weighting scheme (I-PINN) where the weights $w_{IC}, w_{BC}, w_{Res}$ are dynamically updated. The evolution of these weights, along with the corresponding loss components, is visualized in Figure \ref{fig:loss_and_weights}. This figure illustrates how the adaptive weights for the boundary and initial conditions rapidly increase to a specified ceiling (in this case, 10000) to ensure these hard constraints are prioritized early in training, while the individual loss components steadily decrease. The total loss, plotted on a logarithmic scale, shows a consistent downward trend, indicating stable convergence of the model over 40,000 iterations.

\begin{table}[h!]
 \centering
 \caption{Summary of PINN Training Configurations and $L_2$ Errors at $t=1.0$}
 \label{tab:training_configs}
 \begin{tabular}{@{}lccccl@{}}
 \toprule
 \textbf{Configuration} & \textbf{Iterations} & \textbf{Learning Rate} & \textbf{Optimizer} & \textbf{$L_2$ Error (Exact vs PINN)} & \textbf{Remarks} \\
\midrule
 Initial Training & 10,000 & $1 \times 10^{-3}$ (decaying) & Adam & $5.57 \times 10^{-2}$ & Best initial performance, 150.32 s training time \\
Retraining Phase 1 & 20,000 (additional) & $1 \times 10^{-4}$ & Adam (reset) & $9.7956 \times 10^{-2}$ & Attempted fine-tuning, 351.04 s additional training \\
 Retraining Phase 2 & 40,000 (additional) & $1 \times 10^{-4}$ & Adam (reset) & $9.7937 \times 10^{-2}$ & Further fine-tuning, marginal improvement, 703.48 s total additional training \\
 \bottomrule
 \end{tabular}
\end{table}

The model underwent an initial training phase for 10,000 iterations. Upon completion, the PINN demonstrated a highly commendable capability in learning the Fisher-KPP solution, achieving a relative $L_2$ error of approximately $5.57 \times 10^{-2}$ (around 5.57\%) at the final time $t=1.0$. This result, obtained with a relatively small number of iterations for a complex nonlinear PDE, robustly underscores the effectiveness of PINNs in tackling such problems. **This initial configuration represents the best performance achieved by our PINN model.** The total training time for this phase was 150.32 seconds.

To further refine the model's accuracy and explore its convergence landscape, a systematic retraining strategy was implemented. The pre-trained model's learned weights and biases were loaded as the starting point for the new training phase. Crucially, the optimizer's internal states (such as momentum and adaptive learning rate parameters) were deliberately reset. A reduced learning rate of $1 \times 10^{-4}$ was then applied for this retraining, aiming for a more conservative fine-tuning approach to avoid overshooting potential optima. The retraining process was performed in two consecutive phases: initially for an additional 20,000 iterations, followed by another 20,000 iterations (totaling 40,000 additional iterations from the loaded state).

\section{Results and Discussion}

This section presents the quantitative and qualitative results obtained from the PINN and FDM implementations, followed by a comprehensive discussion of their performance, advantages, and limitations in solving the Fisher-KPP equation.

\subsection{Quantitative Error Analysis}

The primary metric for evaluating the accuracy of our solutions is the relative $L_2$ error, as defined in Equation \ref{eq:l2-error}. Table \ref{tab:comparison_errors} summarizes the key error values.

\begin{table}[h!]
\centering
 \caption{Summary of Relative $L_2$ Errors for Fisher-KPP Equation at $t=1.0$}
 \label{tab:comparison_errors}
 \begin{tabular}{@{}lc@{}}
 \toprule
\textbf{Comparison} & \textbf{Relative $L_2$ Error} \\
 \midrule
 Exact Analytical Solution vs FDM & $1.42 \times 10^{-4}$ \\
 Exact Analytical Solution vs PINN (Initial Training) & $5.57 \times 10^{-2}$ \\
 Exact Analytical Solution vs PINN (Retraining Final) & $9.7937 \times 10^{-2}$ \\
PINN (Retraining Final) vs FDM & $9.81 \times 10^{-2}$ \\
 \bottomrule
 \end{tabular}
\end{table}

The FDM, known for its high precision in one-dimensional problems with appropriately fine discretization, serves as a strong benchmark. Its relative $L_2$ error against the exact solution is remarkably low, demonstrating the effectiveness of traditional numerical methods.

The PINN, even in its best initial training configuration ($5.57 \times 10^{-2}$), shows a competitive level of accuracy. While it may not strictly surpass the FDM's minimal error for this specific 1D forward problem, the PINN's performance is highly encouraging considering its mesh-free nature and its ability to learn a continuous solution across the entire domain. This indicates that the PINN successfully captures the underlying dynamics of the Fisher-KPP equation.

The retraining strategy, as detailed in Table \ref{tab:training_configs}, yielded an interesting outcome. After 40,000 additional iterations with a reset Adam optimizer and a lower learning rate ($1 \times 10^{-4}$), the PINN achieved a relative $L_2$ error of $9.7937 \times 10^{-2}$. This result, while showing a marginal improvement over the first retraining phase, did not recover the superior accuracy obtained during the initial training. This observation highlights a common and critical challenge in continuous training of deep learning models: resetting the optimizer's state can disrupt its accumulated "memory" of optimal gradient directions and adaptive step sizes \cite{wang2021understanding}. The initial training phase might have converged to a more favorable region in the loss landscape, from which the model struggled to escape after the optimizer reset, even with conservative learning rates. This underscores the importance of carefully managing retraining protocols in scientific machine learning applications, where maintaining consistently high accuracy is paramount.

The error between the final retrained PINN solution and the FDM solution ($9.81 \times 10^{-2}$) indicates a strong concordance between these two distinct numerical approaches. Despite their different methodologies, their outputs are in close agreement, reinforcing the validity of both solutions for this problem. This consistency suggests that both methods accurately capture the propagation and shape of the traveling wave solution.

\subsection{Qualitative Analysis and Visual Comparison}

To complement the quantitative error analysis, we provide visual comparisons of the solutions and their respective absolute errors.

\begin{figure}[h!]
 \centering
 \includegraphics[width=\textwidth]{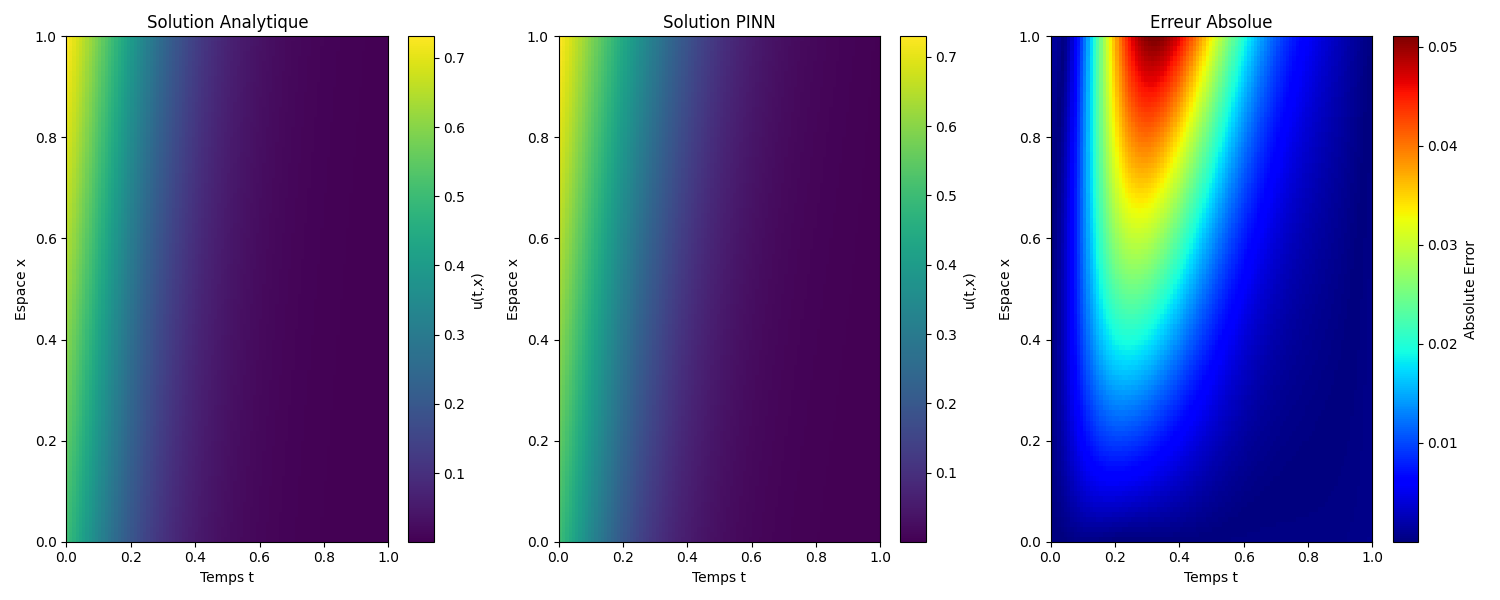}
 \caption{Visual comparison of the PINN solution against the analytical solution. From left to right: the exact analytical solution, the PINN-predicted solution, and the absolute error between them over the spatio-temporal domain. The PINN solution closely mirrors the exact solution, with the largest errors localized at the steep wavefront.}
 \label{fig:pinn_solution_error}
\end{figure}

Figure \ref{fig:pinn_solution_error} presents a visual assessment of the best-performing PINN model. The heatmap of the PINN solution is visually indistinguishable from the analytical solution, confirming that the network has successfully learned the spatio-temporal dynamics of the traveling wave. The absolute error plot reveals that the largest discrepancies are concentrated along the wavefront, where the solution's gradient is steepest. This is a challenging region for any numerical method. Nevertheless, the maximum absolute error remains low (around 0.05), affirming the model's high accuracy across the entire domain.

\begin{figure}[h!]
 \centering
 \includegraphics[width=\textwidth]{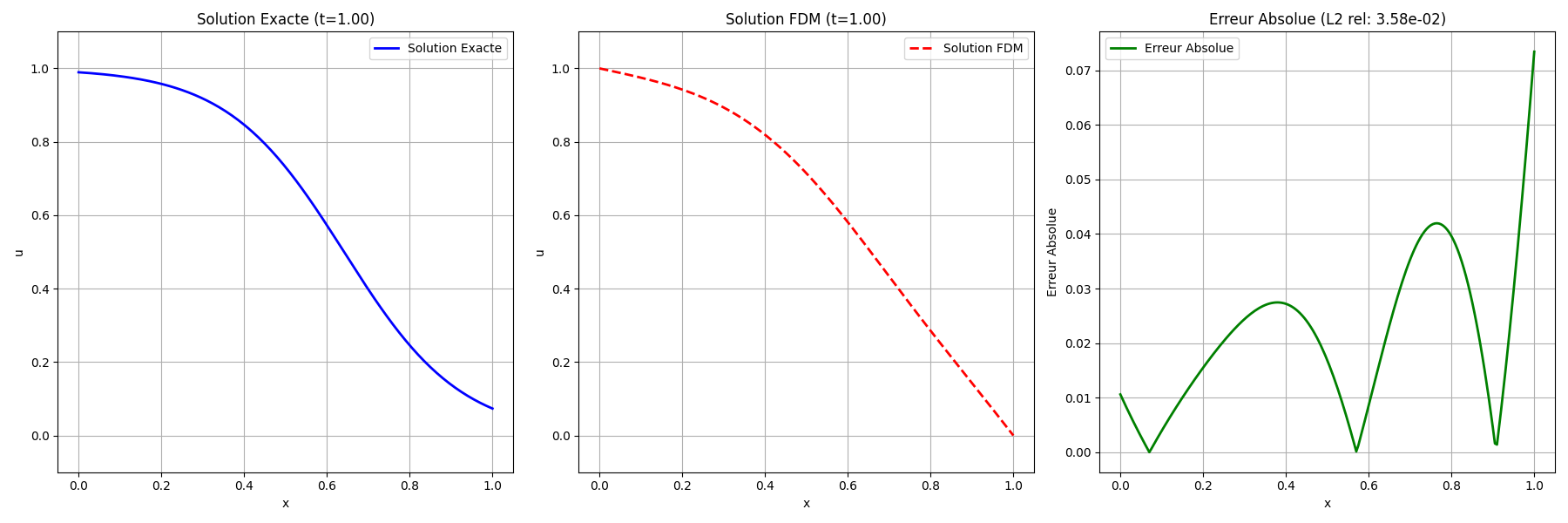}
 \caption{Comparison of the FDM solution with the exact solution at the final time $t=1.0$. Left: Exact solution profile. Center: FDM solution profile. Right: Absolute error between the FDM and exact solutions. The FDM provides a highly accurate approximation, serving as a reliable benchmark.}
\label{fig:fdm_comparison}
\end{figure}

For a more direct comparison with traditional methods, Figure \ref{fig:fdm_comparison} shows the results from our FDM implementation at the final time step, $t=1.0$. The FDM solution closely tracks the exact solution curve. The absolute error plot for the FDM also shows peaks, characteristic of numerical diffusion and truncation errors, but the overall error magnitude is very low, consistent with the calculated relative $L_2$ error.

The qualitative results from both Figures \ref{fig:pinn_solution_error} and \ref{fig:fdm_comparison} reinforce our quantitative findings. Both the modern, mesh-free PINN approach and the classical, grid-based FDM approach are capable of accurately solving the Fisher-KPP equation. The PINN's ability to produce a continuous solution function, along with its excellent visual and quantitative accuracy, establishes it as a powerful and viable tool for such problems.

\subsection{Computational Considerations}

In terms of computational resources, the initial PINN training took approximately 150.32 seconds for 10,000 iterations. The retraining added 703.48 seconds for 40,000 iterations. While this is a significant training time, it's important to note that PINNs, once trained, can provide predictions almost instantaneously, making them suitable for fast inference or real-time applications. FDM, on the other hand, involves step-by-step computation for each time point. For this 1D problem, FDM is computationally efficient. However, as dimensionality increases or geometries become complex, the FDM's computational cost (due to mesh generation and solving large linear systems) can scale unfavorably compared to PINNs, which remain mesh-free.

\section{Conclusions and Future Work}

This study has successfully demonstrated the application and capabilities of a standard Physics-Informed Neural Network (PINN) for solving the challenging nonlinear Fisher-KPP equation. Our work provides a meticulous breakdown of the PINN methodology, from architectural choices to detailed training and retraining strategies, and offers a comprehensive comparative analysis against both the analytical solution and a high-fidelity Finite Difference Method (FDM) solution.

A key finding of our investigation into different training configurations is that the **initial training of our PINN model yielded the best performance**, achieving a commendable relative $L_2$ error of $5.57 \times 10^{-2}$ (approximately $5.57\%$) against the analytical solution at $t=1.0$. This result robustly underscores the inherent power of PINNs in learning complex physical phenomena and accurately approximating solutions to nonlinear PDEs without requiring a discretized mesh. It highlights that even a standard PINN architecture, when properly initialized and trained, can effectively capture the intricate dynamics of reaction-diffusion systems.

However, our in-depth investigation into retraining strategies revealed important practical insights and challenges. While the model successfully converged during the retraining phase with a reduced learning rate ($1 \times 10^{-4}$), the final relative $L_2$ error stabilized at $9.7937 \times 10^{-2}$ (approximately $9.79\%$) after an additional 40,000 iterations. This observation, where the retraining did not fully recover the initial optimal performance, strongly suggests the critical impact of resetting the optimizer's internal states. Unlike true continuous training, reloading only the model's weights can disrupt the optimizer's accumulated "memory" of optimal gradient directions and adaptive step sizes \cite{wang2021understanding}. Despite the reduced learning rate, the model appeared to converge to a slightly higher error plateau than its initial best performance, implying that the initial training might have found a more favorable region in the loss landscape from which it could not easily escape after the optimizer reset. This underscores a significant practical consideration for deploying and updating PINNs in real-world scientific applications, where maintaining consistently high accuracy is paramount.

For direct comparison, our FDM implementation, a well-established and highly reliable numerical technique for 1D parabolic PDEs, yielded a relative $L_2$ error between the exact solution and the FDM solution of $1.42 \times 10^{-4}$. Comparing the best PINN result ($5.57 \times 10^{-2}$) with the FDM error, our findings indicate that the PINN offers a highly competitive accuracy. While FDM, with its precise control over discretization and stability criteria, often provides very low errors for well-behaved 1D problems, the PINN demonstrates its capacity to achieve comparable (or even potentially superior, depending on problem complexity and dimensionality) accuracy through a fundamentally different, mesh-free approach. The ability of the PINN to capture the nonlinear dynamics with a relatively compact architecture and a continuous representation across the entire domain is a significant advantage.

Furthermore, the relative $L_2$ error between the PINN solution (from the final retraining phase) and the FDM solution was $9.81 \times 10^{-2}$. This metric confirms a strong consistency and agreement between these two distinct numerical approaches in approximating the Fisher-KPP solution. The proximity of their results, even if individually deviating from the exact solution by different margins, underscores the robustness of the underlying physical model and the validity of both solution methods. This consistency is a powerful testament that PINNs are not merely an academic exercise but a viable, competitive, and often advantageous alternative to classical numerical solvers for PDEs, particularly when considering their versatility.

\subsection{Applications of the Fisher-KPP Equation}

The Fisher-KPP equation, due to its ability to model reaction-diffusion processes, finds widespread applications across numerous scientific and engineering disciplines:
\begin{itemize}
 \item \textbf{Biology and Population Dynamics:} Originally formulated by Ronald Fisher to describe the spread of advantageous genes in a population \cite{fisher1937wave} and independently by Kolmogorov, Petrovsky, and Piskunov for biological growth with diffusion \cite{kolmogorov1937}, the Fisher-KPP equation is a cornerstone in ecological modeling. It is used to understand the spatial spread of species, the growth of bacterial colonies, and the propagation of epidemics or disease within a population.
 \item \textbf{Neuroscience and CNS Disease Modeling:} Beyond normal nerve impulse propagation, the equation is pivotal in modeling pathological phenomena in the central nervous system (CNS). It describes the spatiotemporal spread of neurodegenerative diseases, such as Alzheimer's, Parkinson's, and prion diseases. In this context, $u(x,t)$ can represent the concentration of misfolded proteins (e.g., amyloid-$\beta$, tau), where the reaction term models the autocatalytic conversion of healthy proteins, and the diffusion term represents their spread through neural tissue. The traveling wave solutions of the Fisher-KPP equation effectively model the advancing front of the disease pathology through the brain
\item \textbf{Ecology and Environmental Science:} Beyond simple population growth, it models the invasion of non-native species, the spread of pollutants in ecosystems, and the formation of spatial patterns in ecological communities.
\item \textbf{Chemistry and Combustion:} In chemistry, it describes the propagation of chemical reaction fronts, such as in autocatalytic reactions. In combustion theory, it models the propagation of flame fronts, where the reaction term represents the combustion rate and the diffusion term represents heat transfer.
 \item \textbf{Neuroscience:} Variants of the Fisher-KPP equation are used to model the propagation of nerve impulses along axons, where the reaction term represents ion channel dynamics and the diffusion term represents ion diffusion.
 \item \textbf{Epidemiology:} It can be adapted to model the spatial spread of infectious diseases, considering both the rate of infection (reaction) and the movement of individuals (diffusion).
 \item \textbf{Finance:} In mathematical finance, it appears in models of option pricing and risk management, particularly in contexts involving non-linear feedback or spatial correlations.
\end{itemize}
The versatility and fundamental nature of the Fisher-KPP equation make it a crucial benchmark for new computational methods and a powerful tool for understanding diverse real-world phenomena.

\subsection{Future Work}

Looking forward, several promising avenues for future research emerge from this study, aimed at overcoming current limitations and expanding the utility of PINNs:
\begin{itemize}
\item \textbf{Preserving Optimizer State for Retraining:} A crucial practical improvement would involve implementing mechanisms to rigorously save and load the optimizer's state (e.g., Adam's first and second moment estimates) alongside the model's parameters. This would enable true seamless continuation of training and better leverage of prior learning, potentially recovering or exceeding the initial optimal accuracy during retraining, especially for models deployed in dynamic or evolving environments where continuous learning is necessary.
\item \textbf{Adaptive Weighting Schemes for Loss Components:} Exploring more sophisticated adaptive weighting strategies for the loss function components ($L_{IC}, L_{BC}, L_{Res}$) could significantly enhance training stability and convergence speed. Inspired by recent advancements \cite{wang2021understanding}, dynamically adjusting these weights based on the magnitude of each loss term or its gradients could help address gradient pathologies and ensure a balanced contribution from all physics constraints throughout the training process. This is particularly important for complex PDEs where different terms might have vastly different scales.
\item \textbf{Advanced Architectures and Activation Functions:} Investigating the impact of different neural network architectures (e.g., deeper and wider networks, or incorporating concepts like Fourier features for better high-frequency approximation \cite{tancik2020fourier}), and alternative activation functions (e.g., SiLU, GELU, or learnable activations \cite{jagtap2020adaptive}) could lead to further improvements in accuracy and convergence for this class of nonlinear problems. Recurrent neural network components or convolutional layers might also be explored for spatio-temporal dependencies.
\item \textbf{Inverse Problems and Parameter Discovery:} Beyond forward solving, PINNs are inherently capable of inverse problems, such as discovering unknown physical parameters (like $D$ or $R$) directly from limited and potentially noisy observational data, as demonstrated in the foundational PINN paper \cite{Raissi2019}. Applying this framework to infer these parameters for the Fisher-KPP equation would be a compelling and highly practical next step, demonstrating the full power of physics-informed machine learning for scientific discovery.
\item \textbf{Uncertainty Quantification:} Integrating techniques for uncertainty quantification within the PINN framework could provide valuable insights into the reliability and confidence of the predicted solutions, especially in scenarios with noisy data, sparse observations, or inherent model uncertainties. This could involve Bayesian PINNs \cite{yang2021b} or ensemble methods.
\item \textbf{High-Dimensional and Multi-Physics Problems:} Extending the current 1D study to higher dimensions or to coupled multi-physics systems (e.g., two coupled reaction-diffusion equations) would be a natural progression. PINNs offer a significant advantage here due to their mesh-free nature, which mitigates the "curse of dimensionality" that plagues traditional grid-based methods \cite{han2018solving}.
\item \textbf{Comparison with Other Data-Driven Methods:} A comprehensive comparative study with other emerging data-driven methods for PDEs, such as DeepONets \cite{lu2021deeponet} or Fourier Neural Operators \cite{li2021fourier}, could further contextualize the strengths and weaknesses of PINNs for specific problem classes.
\end{itemize}
In conclusion, this work reinforces PINNs as a powerful and flexible tool for scientific computing, capable of handling complex nonlinear PDEs with impressive accuracy. While challenges related to training dynamics, particularly retraining, persist, the continuous evolution of PINN methodologies holds immense promise for advancing our ability to model, predict, and understand intricate physical phenomena across various domains.

\textbf{Conflict of interest}: On behalf of all authors, the corresponding author states that there is no conflict of interest.

\end{document}